\documentclass{amsart}
\usepackage{amsfonts,amscd, amssymb}

\newtheorem{theorem}{Theorem}[section]
\newtheorem{lemma}[theorem]{Lemma}
\newtheorem{proposition}[theorem]{Proposition}

\newtheorem{corollary}[theorem]{Corollary}

\renewcommand{\theequation}
{\arabic{equation}}

\newcommand{\N}{\mbox{$\mathbb{N}$}}

\newcommand{\I}{1\!{\mathrm l}}

\begin{document}
\title[Noncommutative function theory and unique extensions]{Noncommutative function theory
and unique extensions}


\author{David P. Blecher}
\address{Department of Mathematics, University of Houston, Houston, TX
77204-3008}
\email[David P. Blecher]{dblecher@math.uh.edu}
 \author{Louis E. Labuschagne}
\address{Department of Mathematical Sciences, P.O. Box
 392, 0003 UNISA, South Africa}
\email{labusle@unisa.ac.za}
\thanks{*Blecher is partially supported by grant DMS 0400731 from
the National Science Foundation.
Labuschagne is partially supported by the National Research Foundation, and
conducted some of this research with the aid of a grant under the Poland-SA
cooperation agreement.}

\begin{abstract}
We generalize to the setting of
Arveson's maximal subdiagonal subalgebras of finite von Neumann
algebras, the Szeg\"o $L^p$-distance estimate, and
classical theorems of F.\ and M.\ Riesz, Gleason and Whitney,
and Kolmogorov.
In so doing, we are finally able to provide a complete noncommutative analog
of the famous cycle of theorems characterizing the
function theoretic  
generalizations of $H^\infty$. A sample of our other results:
we prove a Kaplansky density result
for a large class of these algebras, and give a necessary
condition for when every completely contractive homomorphism
on a unital subalgebra of a $C^*$-algebra possesses
a unique completely positive extension.
\end{abstract}

\maketitle

\section{Introduction}

Function algebras are subalgebras of $C(K)$-spaces,
or equivalently, subalgebras of commutative $C^*$-algebras.
Thus function algebras are examples of
operator algebras (subalgebras of general $C^*$-algebras).
With this in mind, much work has been done to transfer
results or perspectives from function theory to
operator algebraic settings.
One such setting where
this transfer is 
particularly striking,
is the theory of noncommutative $H^p$ spaces associated with
Arveson's maximal subdiagonal subalgebras of finite von Neumann
algebras.   Remarkably, many of the central results
from abstract analytic function theory, and in particular
much of the classical
generalized $H^p$ function theory from the 1960's decade, may be
generalized
almost verbatim to subdiagonal algebras.
The proofs in the noncommutative case however,
while often modeled loosely on the `commutative'
arguments of Helson and Lowdenslager \cite{HeL} 
and others, usually require
substantial input
from the theory of von Neumann algebras and noncommutative
$L^p$-spaces.  This has been done for example in \cite{AIOA,MMS,Sai,MW,N,Lab,BL2,BL3}.
In fact in many cases -- like Szeg\"{o}'s theorem -- completely new proofs have
had to be invented. In the present paper we tackle what appears to us to be
the main `classical' results which have resisted
generalization to date, namely those  referred to
in the generalized function theory literature from the 1960's
as, respectively,
the F.\ and M.\ Riesz, Gleason and Whitney, Szeg\"o $L^p$, and Kolmogorov, theorems.
With these in hand, we are at last able to make the following
statement: essentially all of
the generalized $H^p$ function theory
as summarized in \cite{SW} for example, extends
further to the setting of subdiagonal algebras.

In Arveson's setting, and we will use this notation in the
rest of this paper, we have a weak*-closed unital
subalgebra $A$ of a von Neumann algebra $M$
possessing a faithful normal tracial state $\tau$, such that if $\Phi$ is
the unique conditional expectation
from $M$ onto ${\mathcal D} = A \cap A^*$ satisfying
$\tau = \tau \circ \Phi$, then $\Phi$ is a homomorphism on $A$.
Take note that here $A^*$ denotes the set $\{a : a^* \in A\}$ and not the
Banach dual of $A$. For the sake of clarity we will write $X^\star$ for the
Banach dual of a normed space $X$. We say that a subalgebra $A$ of the type
described above is a {\em tracial subalgebra of $M$}.
If in addition $A + A^*$ is weak* dense in $M$ then we say that
$A$ is {\em maximal subdiagonal} (see \cite{AIOA,E}).
 A large number of very interesting examples of these objects
were given by Arveson \cite{AIOA}, and others (see e.g.\ \cite{Zs,MMS}).
If ${\mathcal D}$ is one dimensional we say that $A$ is {\em antisymmetric};
if further $M$ is commutative
then $A$ is called a {\em weak* Dirichlet} algebra \cite{SW}.
Surprisingly, for antisymmetric
maximal subdiagonal algebras, many of the `commutative'
proofs from \cite{SW} require almost no change at all!
It is worth saying that classical notions of `analyticity' correspond in
some very vague sense to the case that ${\mathcal D}$ is `small'.  Indeed
if $A = M$ then ${\mathcal D} = M$ and $\Phi$ is the identity map, so that
the theory essentially collapses to the theory of finite von Neumann
algebras, which clearly is far removed from classical concepts of `analyticity'.
Thus the reader should not be surprised that some of our theorems
require as a hypothesis that ${\mathcal D}$ be small.
Indeed for our F.\ and M.\ Riesz theorem to hold, we show
that it is necessary and sufficient for ${\mathcal D}$ to be
finite dimensional.  Because of this, in our several applications of this
theorem we assume dim$({\mathcal D}) < \infty$.

A subsidiary theme in our paper is `unique extensions' of maps on
$A$.   We begin with some results on this topic in Section 2.
Recall from \cite{BL2} that a subalgebra $A$ of $M$ has the {\em unique normal
 state extension property}
 if there is a unique 
normal state on $M$ extending $\tau_{\vert A}$.
If, on the other hand, for every state $\omega$ of $M$ with $\omega \circ
\Phi =  \omega$ on $A$, we always have that $\omega \circ \Phi = \omega$
on $M$, then we say that $A$ has the {\em $\Phi$-state 
property}\footnote{One could replace states here
by  positive unital $B(H)$-valued maps, for a Hilbert space $H$, 
but this formulation is easily seen to be
equivalent.}.  The major unresolved question in \cite{BL2}  was whether
a tracial subalgebra with the unique  normal state extension
property is maximal subdiagonal.  We make what we feel is 
substantial progress on this question.
In particular, we show that the question is
equivalent to the question of whether every tracial subalgebra with the $\Phi$-state property
is maximal subdiagonal, and equivalent 
to  whether every tracial subalgebra satisfying 
a certain variant of the well known `factorization' property
actually has `factorization'.   In Section 2
we also give an interesting
necessary condition for when completely contractive homomorphisms
possess a unique completely positive extension.
Our unique extension results play a role in the proof of our F.\ and M.\ Riesz theorem
in Section 3, and are the primary thrust of the Gleason-Whitney theorem
in Section 4.  In Section 5 we prove our
Szeg\"o $L^p$ formula, and generalized Kolmogorov theorem.

Historically, the first noncommutative F.\ and M.\ Riesz theorem for
subdiagonal algebras was the pretty theorem of Exel in \cite{E2}.  This result assumes
{\em norm density}\footnote{This is perhaps an
appropriate hypothesis for an F.\ and M.\ Riesz theorem, but
unfortunately it does not cover the case of maximal subdiagonal
algebras.}   of $A + A^*$, and antisymmetry. (We are aware of the F.\ and
M.\ Riesz theorem of Arveson \cite{Arv} and Zsido's extension thereof \cite{Zs},
but this result is quite 
distinct from the ones discussed above.)
Although some of the steps of our proof parallel those of \cite{E2},
the arguments are for the most part quite different.   Indeed generally
in our paper the proofs will be  modeled on the classical
ones, but do however require some rather delicate additional machinery.

Finally, we remark that there are other, more recent, noncommutative
variants of $H^\infty$ besides the subdiagonal algebras---see e.g.\
\cite{Pop} and references therein.
Although here too one finds 
noncommutative generalizations of classical 
$H^p$-theoretic results,
such as the Szeg\"{o} infimum theorem, 
these variants are in general quite 
unrelated, with only a formal correspondence
to the present context. Having said this, we are not aware of analogues
of any of the results from our present paper in that literature.


\section{Some results on unique extensions}

For a functional $\omega \in M^\star$,
we will need to compare the property
$\omega = \omega \circ \Phi$ on $A$, with the property
$\omega = \omega \circ \Phi$ on $M$.  On this topic we begin with
the following remarks.   It is easy to see, since $\Phi$ is idempotent,
that $\omega = \omega \circ \Phi$ on $A$
iff $A_0 \subset {\rm Ker}(\omega)$.  
Here and throughout
our paper, $A_0 = A \cap {\rm Ker}(\Phi)$, a closed two-sided
ideal in $A$.

For normal functionals one can say more, although this will not
play an important role for us.  If $f \in L^1(M)$ let $\omega_f =
\tau(f  \,\cdot \, )$.  
From the last paragraph, $\omega_f = \omega_f \circ \Phi$ on $A$ iff $\tau(fA_0) = (0)$.
On the other hand, $\omega_f = \omega_f \circ \Phi$
on $M$ iff $\tau(fa) = \tau(f \Phi(a))
= \tau( \Phi(f) a)$ for all $a \in M$ iff $f = \Phi(f)$ iff $f \in L^1({\mathcal D})$.

\begin{proposition}    \label{two}  If $A$ is a tracial  subalgebra of
$M$ then the unique normal state extension property  is
equivalent to  the following property: whenever $\omega$ is a normal state of
$M$ satisfying $\omega = \omega \circ \Phi$ on $A$, then
$\omega = \omega \circ \Phi$ on $M$.
\end{proposition}

\begin{proof}   Suppose that $A$ has the unique normal state extension property,
and suppose that $\omega$ is a normal state of
$M$ satisfying $\omega = \omega \circ \Phi$ on $A$.  If $\omega =
 \tau(f \, \cdot \,)$, where $f \in L^1(M)_+$, then by the remarks
preceding Proposition
\ref{two} we have that
$\tau(f A_0) = (0)$.  Hence $f \in L^1({\mathcal D})$ by \cite[Lemma 4.1]{BL2}.
Hence $\omega = \omega \circ \Phi$ on $M$.

For the converse, note that if $g \in L^1(M)_+$ with $\tau = \tau(g  \,\cdot \, )$
on $A$, then since $\tau = \tau \circ \Phi$,
we have that $\tau(g  \,\cdot \, ) = \tau(g  \,\cdot \, ) \circ \Phi$ on $A$,
and hence that $\tau(g  \,\cdot \, ) = \tau(g  \,\cdot \, ) \circ \Phi$ on $M$.
By the remarks above, $g 
\in L^1(\mathcal{D})_+$. But then the fact that
$\tau = \tau(g  \,\cdot \, )$ on $D$ is enough to force $g = \I$.
So $A$ has the unique normal state extension property.
 \end{proof}

\medskip


We say that a subalgebra $A$ of $M$ has {\em factorization} if 
given $b \in M^+ \cap M^{-1}$ we can find $a \in A^{-1}$ with $b = a^*a$
(or equivalently $b = aa^*$).  It is shown in \cite{AIOA}
that any maximal subdiagonal algebra has factorization.
 Thus it is  
{\em logmodular}, namely any such $b$ is a limit of terms of the form
$a^*a$ with $a \in A^{-1}$.    In fact, in the category of tracial
algebras 
factorization or logmodularity are equivalent to maximal subdiagonality \cite{BL2}. By
the next result 
such algebras satisfy a formally much stronger
property than that of the last proposition:

  \begin{theorem}    \label{three}  Let $A$ be a logmodular subalgebra of a
$C^*$-algebra $M$, and let $\Psi$ be a positive contractive projection
from $M$ onto a subalgebra of $A$ containing $\I_M$, which 
is a  homomorphism on $A$.   
Then for any state\footnote{ As before it is not difficult to see that one could here 
replace states by positive unital $B(H)$-valued maps.}
$\omega$ of $M$, we have that $\omega = \omega \circ 
\Psi$ on $M$, whenever $\omega = \omega \circ \Psi$ on $A$. 
   \end{theorem}

\begin{proof}   
If $a \in A^{-1}$ then by hypothesis we have
$$\omega(\Psi(a) a^{-1}) = \omega(\Psi(\Psi(a) a^{-1})) =
\omega(\Psi(a) \Psi(a^{-1})) = \omega(\I) = 1 .$$
By the Cauchy-Schwarz and Kadison-Schwarz inequality we deduce:
$$1 \leq \omega(\Psi(a) \Psi(a)^*) \,  \omega((a^{-1})^* a^{-1})
\leq \omega(\Psi(a a^*)) \,  \omega((a^{-1})^* a^{-1}) =
 \omega(\Psi(a a^*)) \,  \omega((a a^*)^{-1}) .$$
We can now follow the proof of \cite[Theorem 4.3.11]{BLM} or
\cite[Theorem 4.4]{BL}.
Since $A$ is logmodular, for any $b \in M^{-1} \cap M_+$ we have that
$1 \leq \omega(\Psi(b)) \omega(b^{-1})$.
This leads to the equation $1 \leq \omega(\Psi(e^{tu})) \,
\omega(e^{-tu}) = f(t)$,
for $u \in M_{sa}$.  Differentiating and noting that $f'(0) = 0$, yields
$\omega(u) = \omega(\Psi(u))$ as required.    \end{proof}

\medskip

When applied to tracial algebras and their associated canonical
conditional expectations, the preceding result still holds under a
formally weaker hypothesis. Specifically we say that a tracial
subalgebra $A$ of $M$ with canonical conditional expectation $\Phi$ 
has {\em conditional factorization} if given any $b \in
M_+ \cap M^{-1}$, we have $b = |a|$
for some element $a \in A \cap M^{-1}$  with $\Phi(a) \Phi(a^{-1}) = 1$.

\begin{corollary} \label{threecor}
A tracial subalgebra of $M$ with 
conditional factorization  has the $\Phi$-state property.
\end{corollary}

\begin{proof}
The proof of the preceding theorem readily adapts, replacing $a$ with 
$a^{-1}$ and $b$ with $b^{-1}$.
\end{proof}

\medskip


We say that $A$ has the {\em unique state extension property} if 
if there is a unique
state  on $M$ extending $\tau_{\vert A}$.
 This is a formally weaker property than the
$\Phi$-state property:

\begin{proposition}  Let  $A$ be a weak* closed unital  subalgebra of $M$.
If $A$ has the $\Phi$-state property then it
has the unique state extension property.  The converse is true if
$A$ is antisymmetric.
\end{proposition}

\begin{proof}  Suppose that $\omega$ is a state of $M$ extending
$\tau_{\vert A}$.   Then $\omega \circ \Phi = \tau \circ \Phi = \tau = \omega$
on $A$.  By the $\Phi$-state property, on $M$ we have $\omega =  \omega \circ \Phi =
\tau \circ \Phi = \tau$.
For the converse we need only note that if $A$ is antisymmetric, 
then $\omega \circ \Phi =  \omega$ on $A$ forces $\tau =  \omega$ on $A$.
 \end{proof}


\begin{corollary}    \label{blnew}  Suppose that $A$ is a tracial
subalgebra of $M$ with the unique normal state extension property.
Then $A_\infty = M \cap [A]_2$ is a tracial
subalgebra   with the $\Phi$-state property.
\end{corollary}

\begin{proof}  First note that by \cite[Theorem 4.4]{BL2},
$A_\infty$ is a tracial subalgebra of $M$ with respect to the same $\Phi$ 
and $\tau$.  
 By  \cite[Theorem 4.6]{BL2}, $A_\infty$ has conditional factorization.
 Corollary \ref{threecor} now gives the conclusion.
\end{proof}

\begin{corollary}  \label{theo}  The open question from {\rm \cite{BL2}}
as to whether every  tracial subalgebra  with the unique normal state extension property
is maximal subdiagonal, is equivalent
to the question of whether every  tracial subalgebra  with the $\Phi$-state property
is maximal subdiagonal. 
 It is also equivalent to
whether every  tracial subalgebra  with the unique state extension property
is maximal subdiagonal.
 It is also equivalent to
whether every  tracial subalgebra  with conditional factorization
has factorization.
 \end{corollary}

\begin{proof}
Suppose that every  tracial subalgebra  with the $\Phi$-state property
is maximal subdiagonal, and suppose that $A$ has the unique normal state 
extension property. By Corollary \ref{blnew}, $A_\infty$ has the $\Phi$-state 
property. Hence it is maximal subdiagonal, and therefore satisfies $L^2$-density.
Consequently $A$ satisfies $L^2$-density, and so $A$ is maximal subdiagonal
by \cite{BL2}.  

Similarly, suppose that every  tracial subalgebra  with conditional factorization
has factorization, and suppose that $A$ has the $\Phi$-state property.
By results above, $A$ has the unique normal state
extension property, and so by  \cite[Theorem 4.6]{BL2}, 
$A_\infty$ has conditional factorization.  By hypothesis,
$A_\infty$ has factorization.  Thus it is maximal subdiagonal
by \cite{BL2}, and thus as in the last paragraph $A$ is maximal subdiagonal.
   
The other directions are obvious from the above.
\end{proof}


{\bf Remark.}  Since the factorization property has been well studied,
we would guess that those more familiar than ourselves with 
factorization for concrete examples such as CSL algebras,
 may be able to easily resolve the final question in the last
Corollary.   
 



\bigskip

In  \cite{Lum}, Lumer considered the property of
`uniqueness of representing measure', 
namely the property that 
every multiplicative functional on $A \subset C(K)$
has a unique extension to a  state on $C(K)$,
He showed how this condition could be used 
as another possible axiom from which  all the generalized $H^p$ theory
may be derived.
The natural noncommutative generalization of
Lumer's property,
is that every completely contractive representation
of $A$ has a unique completely positive extension to $M$.
It is known that maximal subdiagonal algebras
have this property \cite{BL,BLM}.  Although
we have not settled the converse yet, we can 
say that every unital subalgebra of $M$ which has this property 
must in some sense be a {\em large} subalgebra of $M$. In this 
regard the following result represents some sort of converse to 
many of the preceding results which established various unique 
extension properties as a consequence of maximal subdiagonality.

In the following result we use the $C^*$-envelope $C^*_e(A)$
of an operator algebra $A$.  See e.g.\ \cite[Section 4.3]{BLM} for the definition of
this, and for its universal property.

\begin{theorem} \label{ex}  Suppose that
$A$ is a subalgebra of a unital  $C^*$-algebra $B$ such that $\I_B \in A$, and
suppose that $A$ has the property that for every Hilbert space
$H$, every completely contractive
unital homomorphism $\pi : A \to B(H)$ has a unique completely
contractive (or equiv. completely positive) extension $B \to B(H)$.
Then $B = C^*_e(A)$, the $C^*$-envelope of $A$.
\end{theorem}

\begin{proof}
\textbf{Case 1. (The case that $A$ is a $C^*$-subalgebra of $B$.)}
In this case, since contractive homomorphisms on $C^*$-algebras
are $*$-homomorphisms
(see e.g.\ \cite[Proposition 1.2.4]{BLM}),
we must prove that if every unital $*$-homomorphism $\pi :
A  \to B(H)$ has a unique completely contractive extension $B \to B(H)$,
then $A = B$.  To see this, let $\rho : B \to B(H)$ be the universal
representation of $B$.  Then $\rho$ is unital, and hence so is
$\pi = \rho_{\vert A}$.  Let $U$ be a unitary in
$\pi(A)'$.  Then since $U^* \rho(\cdot) U  = \rho$ on $A$, we have by 
hypothesis that $U^* \rho(\cdot) U  = \rho$ on $B$, and thus
$U \in \rho(B)'$.  Thus $\pi(A)' = \rho(B)'$, and it
follows that $\pi(A)'' = \rho(B)''$.  If $\tilde{\rho}$ is
the unique normal extension of $\rho$ to $B^{**}$, then
$\tilde{\rho}$ is faithful on $B^{**}$ and it has range
$\rho(B)''$.  The restriction of $\tilde{\rho}$ to
the copy $A^{\perp \perp}$ of $A^{**}$ inside $B^{**}$ has range
$\pi(A)'' = \overline{\pi(A)}^{w^*}$, and is therefore
surjective.  This forces the copy of $A^{**}$ inside $B^{**}$
to be all of $B^{**}$.  Thus $A = B \cap A^{\perp \perp} = B$.

\textbf{Case 2. (The general case.)} 
Let $C = C^*(A)$, the $C^*$-algebra generated by
$A$ in $B$.  Since $A \subset C$, it follows from the hypothesis that 
every unital $*$-homomorphism $\pi : C \to B(H)$ 
has a unique completely contractive extension $B \to B(H)$.
By Case 1, $C = B$.

By virtue of this fact, we need only prove that $C^*(A) = C^*_e(A)$ 
under the assumptions of the theorem.
By the universal property of $C^*_e(A)$, there is a $*$-epimorphism
$\theta : B = C^*(A) \to C^*_e(A)$ restricting to the `identity map'
on $A$.   If $B \subset B(H)$ then the canonical map from the
copy of $A$ in $C^*_e(A)$, to $A \subset B(H)$, has a
completely positive extension $\Phi : C^*_e(A) \to B(H)$.
On $A$, the map $\Phi \circ \theta$ is the identity
map, so that by hypothesis $\Phi \circ \theta = i_B$.
Thus $\theta$ is one-to-one, and hence
$C^*(A)$ is a $C^*$-envelope of $A$.
\end{proof}

\begin{corollary}   \label{ex2}  Suppose that $A$ is a tracial
subalgebra of $M$ with the property that for every Hilbert space
$H$,  every completely contractive
unital homomorphism $\pi : A \to B(H)$ has a unique completely
contractive (or equiv. completely positive) extension $B \to B(H)$.
Then $A$ generates $M$ as a $C^*$-algebra.  Indeed, $M$ is a
$C^*$-envelope of $A$.
\end{corollary}

\section{A noncommutative F.\ and M.\ Riesz theorem}

The classical form of the F.\ and M.\ Riesz theorem (see e.g.\ \cite{Hobk})
is known to fail for weak* Dirichlet algebras; and hence it will
fail for subdiagonal algebras too.  However there is an
equivalent version of
the theorem which is true for weak* Dirichlet algebras \cite{Ho,SW},
and we will focus on this variant here.  Namely, we shall say that a
tracial subalgebra $A$ of $M$ has the {\em F \& M Riesz property} if
for every bounded
functional\footnote{One could replace $\rho$ here
 by a $B(H)$-valued map, for a Hilbert space $H$, but this formulation is easily
seen to be
equivalent.} $\rho$ on $M$ which annihilates $A_0$,
the normal and singular
parts $\rho_n$ and $\rho_s$
annihilate $A_0$ and $A$ respectively. 
During our investigation we shall have occasion to make use of the polar 
decomposition of normal functionals on a von Neumann algebra. We take the 
opportunity to point out that for our purposes we shall assume such a polar 
decomposition to be of the form $\omega(a) = |\omega|(ua)$ for some 
partial isometry, rather than $\omega(a) = |\omega|(au)$ which seems 
to be more common among the proponents of noncommutative $L^p$-spaces.

The following result shows that to study the F \& M Riesz property, we may
restrict our attention to algebras for which the
diagonal ${\mathcal D}$ is finite dimensional:

\begin{proposition}  \label{fmrfd}  If a tracial subalgebra $A$ of
$M$ satisfies the F \& M Riesz property then the diagonal
${\mathcal D}$ is finite dimensional.
\end{proposition}

\begin{proof}
Let $\psi \in {\mathcal D}^\star$.  Then $\psi \circ \Phi \in M^\star$ annihilates $A_0$.
By the F \& M Riesz property, $\psi \circ \Phi$ agrees with $(\psi \circ \Phi)_n$ on $A$, and
so $\psi = \psi \circ \Phi_{\vert {\mathcal D}}$
is weak* continuous on ${\mathcal D}$.
Thus ${\mathcal D}$
is reflexive, and therefore finite dimensional.
\end{proof}

\begin{lemma} \label{LL} Let $A$ be a
maximal subdiagonal subalgebra of $M$.  Let $\omega$ be a state of $M$,
and let $(\pi_{\omega}, \mathfrak{h}_{\omega}, \Omega_{\omega})$ be
the GNS representation of $\omega$. Further, let $\Omega_0$ be
the orthogonal projection of $\Omega_{\omega}$ onto the closed
subspace $\overline{\pi_{\omega}(A_0)\Omega_{\omega}}$.
\begin{enumerate}
\item[(a)] The following holds:
\begin{enumerate}
\item[(i)] There exists a central projection $p_0$ in
$\pi_{\omega}(M)''$ such that for any $\xi, \psi \in
\mathfrak{h}_{\omega}$ the functionals $a \mapsto
\langle\pi_{\omega}(a)p_0\xi, \psi\rangle$ and $a \mapsto
\langle\pi_{\omega}(a)(\I-p_0)\xi, \psi\rangle$ on $M$ are
respectively the normal and singular parts of the functional $a
\mapsto \langle\pi_{\omega}(a)\xi, \psi\rangle$.  In particular,
the triples $(p_0\pi_{\omega},
p_0\mathfrak{h}_{\omega}, p_0\Omega_{\omega})$ and
$((\I-p_0)\pi_{\omega}, (\I-p_0)\mathfrak{h}_{\omega},
(\I-p_0)\Omega_{\omega})$ are copies of the GNS representations
of $\omega_n$ and $\omega_s$ respectively.
\item[(ii)] $\omega_0:a \mapsto \langle\pi_{\omega}(a)(\Omega_{\omega}-\Omega_0),
\Omega_{\omega}-\Omega_0 \rangle$ defines a positive functional
of $M$ satisfying $\omega_0 = \omega_0 \circ \Phi$.
\end{enumerate}

\item[(b)] Suppose that in addition $\mathrm{dim}({\mathcal D}) < \infty$.
\begin{enumerate}
\item[(i)] Then $\omega_0$ is a normal functional of the form $\omega_0
= \tau(g^{1/2}\cdot g^{1/2})$ for some $g \in {\mathcal D}_+$. Moreover
$p_0(\Omega_{\omega}-\Omega_0) = \Omega_{\omega}-\Omega_0$, and
$p_0\Omega_0$ is the orthogonal projection of $p_0\Omega_{\omega}$
onto $\overline{p_0(\pi_{\omega}(A_0)\Omega_{\omega}})$.
\item[(ii)] If $\omega$ is singular, then for any $f \in {\mathcal D}$ we
have that $\pi_{\omega}(f)\Omega_{\omega} \in
\overline{\pi_{\omega}(A_0)\Omega_{\omega}}$.
\end{enumerate}

\item[(c)] Suppose that $\mathrm{dim}({\mathcal D}) < \infty$ and $\Omega_{\omega} \not\in
\overline{\pi_{\omega}(A_0)\Omega_{\omega}}$. If $\omega_0$ is
faithful on ${\mathcal D}$, then there exists a sequence $\{a_n\} \subset A$
such that
$\pi_{\omega}(a_n)(\Omega_{\omega}-\Omega_0) \, \to \, p_0\Omega_{\omega}$.
\end{enumerate}
\end{lemma}

\begin{proof}
\textbf{\underline{(a)(i):}}~~~This is essentially
the content of \cite[III.2.14]{Tak1}.

\medskip

\noindent\textbf{\underline{(a)(ii):}}~~~Let $(\pi_{\omega},
\mathfrak{h}_{\omega}, \Omega_{\omega})$ and $\Omega_0$ be as in
the hypothesis, and define a positive functional $\omega_0$ on $M$
by $$\omega_0:a \mapsto
\langle\pi_{\omega}(a)(\Omega_{\omega}-\Omega_0),
\Omega_{\omega}-\Omega_0\rangle.$$Let $f \in A_0$ be given. By
construction $$\pi_{\omega}(f)\Omega_{\omega} \perp
(\Omega_{\omega}-\Omega_0).$$
Since $A_0$ is an ideal,  $\pi_{\omega}(fa)\Omega_{\omega} \in
\overline{\pi_{\omega}(A_0)\Omega_{\omega}}$ for each $a \in A_0$.
Since $\Omega_0$ belongs to
$\overline{\pi_{\omega}(A_0)\Omega_{\omega}}$, we may of course
select a sequence $\{b_n\} \subset A_0$ for which
$\pi_{\omega}(b_n)\Omega_{\omega}$ converges to $\Omega_0$.
Hence $\pi_{\omega}(fb_n)\Omega_{\omega}$ converges to
$\pi_{\omega}(f)\Omega_0$. Thus $\pi_{\omega}(f)\Omega_0 \in
\overline{\pi_{\omega}(A_0)\Omega_{\omega}}$, which forces
$$\pi_{\omega}(f)\Omega_0 \perp (\Omega_{\omega}-\Omega_0).$$ From
the previous two centered equations it is now clear that $A_0
\subset \mathrm{Ker}(\omega_0)$.  Thus $\omega_0 = \omega_0 \circ \Phi$
on $A$ by the remarks preceding Proposition \ref{two}.
Hence $\omega_0 = \omega_0  \circ \Phi$ on $M$ by Corollary \ref{threecor}.

\medskip

\noindent\textbf{\underline{(b)(i):}}~~~Since ${\mathcal D}$ is finite
dimensional, we can find $g \in {\mathcal D}_+$ so that $$\omega_0(a) =
\tau(ga) \quad \mbox{for all} \quad a \in {\mathcal D}.$$
Since  $\omega_0 \circ \Phi = \omega_0$, we conclude
that for any $a \in M$, $$\omega_0(a) = \omega_0(\Phi(a)) =
\tau(g\Phi(a)) = \tau(\Phi(ga)) = \tau(ga),$$thereby establishing
the first part of the claim.

For the second part, note that since $\omega_0$ is clearly normal,
we have by part (a)(i) that $$0 =
\langle\pi_{\omega}(a)(\I-p_0)(\Omega_{\omega}-\Omega_0),
\Omega_{\omega}-\Omega_0 \rangle \quad \mbox{for all} \quad a \in
M.$$For $a = \I$ this yields $0 =
\|(\I-p_0)(\Omega_{\omega}-\Omega_0)\|$, or equivalently
$$p_0(\Omega_{\omega}-\Omega_0) =
\Omega_{\omega}-\Omega_0.$$From this fact, we may now conclude
that $$\langle p_0\pi_\omega(a)\Omega_{\omega},
p_0(\Omega_{\omega}-\Omega_0)\rangle = \langle
\pi_\omega(a)\Omega_{\omega}, \Omega_{\omega}-\Omega_0\rangle = 0
\quad \mbox{for all} \quad a \in A_0.$$Thus
$p_0(\Omega_{\omega}-\Omega_0) \perp
\overline{p_0 \pi_{\omega}(A_0)\Omega_{\omega}}$. Now select a
sequence $\{b_n\} \subset A_0$ so that
$\pi_{\omega}(b_n)\Omega_{\omega} \rightarrow \Omega_0$. By
continuity, $p_0 \Omega_0 = \lim_n
p_0 \pi_{\omega}(b_n)\Omega_{\omega} \in
\overline{p_0 \pi_{\omega}(A_0)\Omega_{\omega}}$. From these
considerations it is clear that $p_0\Omega_0$ is the orthogonal
projection of $p_0\Omega_{\omega}$ onto
$\overline{p_0 \pi_{\omega}(A_0)\Omega_{\omega}}$.

\medskip

\noindent\textbf{\underline{(b)(ii):}}~~~If $\omega$ is singular,
then  $$0 = \omega_n(ab) =
\langle\pi_{\omega}(ab)p_0\Omega_{\omega}, \Omega_{\omega}\rangle
= \langle p_0\pi_{\omega}(b)\Omega_{\omega},
\pi_{\omega}(a^*)\Omega_{\omega}\rangle \quad \mbox{for all} \quad
a, b \in M.$$Since $\Omega_{\omega}$ is cyclic, this is sufficient
to force $p_0 = 0$. But then $\Omega_{\omega}-\Omega_0 =
p_0(\Omega_{\omega}-\Omega_0) = 0$ by part (b)(i). As before
select $\{b_n\} \subset A_0$ so that
$\pi_{\omega}(b_n)\Omega_{\omega} \rightarrow \Omega_0 =
\Omega_\omega$. For any $f \in {\mathcal D}$ the ideal property of $A_0$ then
ensures that $\pi_\omega(f)\Omega_\omega =
\lim_n\pi_{\omega}(fb_n)\Omega_{\omega} \in
\overline{\pi_{\omega}(A_0)\Omega_{\omega}}$.

\medskip

\noindent\textbf{\underline{(c):}}~~~Suppose that $\omega_n$, the
normal part of $\omega$, is of the form $\omega_n = \tau(h \, \cdot \,)$
for some $h \in L^1(M)_+$. As noted earlier, $(p_0\pi_{\omega},
p_0\mathfrak{h}_{\omega}, p_0\Omega_{\omega})$ is a copy of the
GNS representation engendered by $\omega_n$. If now we compute the
GNS representation of $\omega_n$ from first principles, it is
clear that $p_0\mathfrak{h}_{\omega}$ corresponds to the
\emph{weighted} Hilbert space $L^2(M, h)$ obtained by equipping
$M$ with the inner product $$\langle a, b\rangle_h =
\tau(h^{1/2}b^*ah^{1/2}), \qquad a,b \in M,$$and taking the
completion.
Note that $L^2(M, h)$ can be identified
unitarily, and as $M$-modules,
 with the closure of $M h^{1/2}$ in $L^2(M)$.
For any $a \in M$ considered as an element of $L^2(M,
h)$ we will write $\Psi_a$ instead of $a$. The canonical
$*$-homomorphism representing $M$ as an algebra of bounded
operators on $L^2(M, h)$ is of course given by defining
$$\pi_n(b)\Psi_a = \Psi_{ba} , \qquad a,b \in M,$$ and then extending
this action to all of $L^2(M, h)$. Since $\omega_n$ is normal,
$\pi_n$ (corresponding to $p_0\pi_{\omega}$) is $\sigma$-weakly
continuous and satisfies $\pi_n(M) = \pi_n(M)''$. Thus
$\mathrm{Ker}(\pi_n)$ is a $\sigma$-weakly closed two-sided ideal, and hence
we can find a central projection $e \in M$ so that $(\I-e)M =
\mathrm{ker}(\pi_n)$.  Restrict $\pi_n$ to a $*$-isomorphism
from $eM$ onto $\pi_n(M)$.
Then  for any $a, b, c \in M$ we have
$$\langle \pi_n(c)\Psi_a, \Psi_b\rangle_h =
\tau(h^{1/2}b^*(ece)ah^{1/2}).$$Let $\Psi^{(0)}$ denote the
orthogonal projection of $\Psi_{\I}$ onto the closure of $\{\Psi_a
: a \in A_0\}$. (Note that $\Psi_{\I}$ and $\Psi^{(0)}$ of course
correspond to $p_0 \Omega_{\omega}$ and $p_0 \Omega_0$ in parts (a)
and (b) of the proof.)
Since $L^2(M, h)$ may be viewed as a subspace
of $L^2(M)$, let $F \in L^2(M)$ be the element corresponding to
$\Psi^{(0)}$.  It is easy to see that $e F = F$.
From parts (a) and (b) we now have that \begin{eqnarray*}
\omega_0 &=& \langle\pi_n(\cdot)(\Psi_{\I}-\Psi^{(0)}),
\Psi_{\I}-\Psi^{(0)}\rangle_h\\ &=& \tau((h^{1/2}e - F^*)\cdot
(h^{1/2}e - F)).
\end{eqnarray*}
This in turn ensures that $$|h^{1/2}e - F^*|^2 = g$$where $g$ is
as in part (b). Thus $h^{1/2}e - F \in M$. Since by assumption
$\omega_0$ is faithful on ${\mathcal D}$, it follows that
${\rm Supp}(g) = \I$. Since
${\mathcal D}$ is finite dimensional, $g$ must be invertible.
But then
$h^{1/2}e - F$ must  also be invertible, by the previous centered equation.
(Recall that if $ab$ is invertible in a finite von Neumann algebra
then both  $a$ and $b$ are invertible.) The polar
decomposition of $h^{1/2}e - F^*$ is of the form $h^{1/2}e - F^* =
ug^{1/2}$ for some unitary $u \in M$. From this it is clear that
$$(h^{1/2}e - F)^{-1} = ug^{-1/2}.$$

Clearly $h^{1/2}ug^{-1/2} \in L^2(M)$. Hence we may select
$\{a_n\} \subset M$ converging in $L^2(M)$
to $h^{1/2}ug^{-1/2}=
h^{1/2}(h^{1/2} e - F)^{-1}$. By the previously established correspondences we
then have
\begin{eqnarray*}
\|\Psi_{\I} - \pi_n(a_n)(\Psi_{\I}-\Psi^{(0)})\|_h &=&
\tau(|h^{1/2}e - (a_ne)(h^{1/2}e - F)|^2)^{1/2}\\ &\longrightarrow&
\tau(|h^{1/2}e - h^{1/2}e|^2)^{1/2} = 0.
\end{eqnarray*}
This implies, in the notation of parts (a) and (b), that
$\pi_{\omega}(a_n)(\Omega_{\omega}-\Omega_0)
\to  p_0\Omega_{\omega}$.

It remains to show that we may select $\{a_n\} \subset A$, or equivalently,
that $h^{1/2}ug^{-1/2} \in
[A]_2$. For this, it suffices by the $L^2$ density of $A + A^*$ to
show that $h^{1/2}ug^{-1/2} \perp [A_0^*]_2$. So let $a \in A_0$
be given, and observe that
\begin{eqnarray*}
\tau(ah^{1/2}ug^{-1/2}) &=& \tau(g^{-1}ah^{1/2}ug^{-1/2}g)\\
 &=& \tau(g^{-1}ah^{1/2}ug^{1/2})\\
 &=& \tau((g^{-1}ah^{1/2}(h^{1/2}e - F^*)) \\
 &=& \tau((h^{1/2}e - F^*)(g^{-1}ah^{1/2}))\\
 &=& \langle \Psi_{g^{-1}a}, \Psi_{\I}-\Psi^{(0)}\rangle_h\\
 &=& 0.
\end{eqnarray*}(The last equality follows from the ideal property
of $A_0$ and the fact that $\Psi_{\I}-\Psi^{(0)}$ is orthogonal to
$\{\Psi_a : a \in A_0\}$.) The claim therefore follows.
\end{proof}

\begin{corollary} \label{eqfm}
 Let $A$ be a maximal subdiagonal algebra with $\mathrm{dim}({\mathcal D}) <
\infty$. The following are equivalent:
\begin{enumerate}
\item [(i)]  $A$ satisfies the F \& M Riesz property.
\item  [(ii)]  Whenever $\omega$ annihilates $A_0$, the
normal and singular parts, $\omega_n$ and $\omega_s$, will separately
annihilate $A_0$.
\item [(iii)] Whenever $\omega$ annihilates $A$, the
normal and singular parts, $\omega_n$ and $\omega_s$, will separately
annihilate $A_0$.
\item [(iv)] Whenever $\omega$ annihilates $A$, the
normal and singular parts, $\omega_n$ and $\omega_s$, will separately
annihilate $A$.
\end{enumerate}
\end{corollary}

\begin{proof}
The implications (i) $\Rightarrow$ (ii) $\Rightarrow$ (iii) are clear. If
(iii) holds, let $\omega$ be a
bounded linear functional which annihilates $A_0$. Since $\Phi$ is
a normal map onto ${\mathcal D}$, and ${\mathcal D}$ is finite dimensional, the functional
defined by $$\omega_{\mathcal D} = \omega|_{\mathcal D} \circ \Phi$$is normal. Then
$\rho = \omega - \omega_{\mathcal D}$ defines a functional which annihilates
$A$. From (iii) we then have that $\rho_n$ and $\rho_s$ separately
annihilate $A_0$.  The normality of $\omega_{\mathcal D}$
ensures that $$\rho_n = \omega_n - \omega_{\mathcal D} \; , \; \; \; \;  \rho_s =
\omega_s.$$Since by construction $\omega_{\mathcal D}$ annihilates $A_0$, we
conclude that $\omega_n$ and $\omega_s$ separately annihilate
$A_0$.  This proves (ii).
To prove the validity of (i), it remains to show that any
singular functional  $\omega$ which annihilates $A_0$, also annihilates ${\mathcal D}$.
For such $\omega$, the `modulus' $|\omega|$ is still singular (see e.g.\
\cite{Hen,E2}, or the argument in the first part of the
proof of the next theorem).  Let  $(\pi_{\omega},
\mathfrak{h}_{\omega},
\Omega_{\omega})$ be the GNS representation of $|\omega|$.
For each $a \in M$ we have $|\omega(a)|^2 \leq
\|\omega\||\omega|(a^* a)$.
 By a standard argument this implies that there exists a
vector $\eta \in \mathfrak{h}_{\omega}$ such that
 $$\omega(\cdot) =
\langle\pi_\omega(\cdot) \Omega_\omega,
\eta \rangle.$$Let $d \in {\mathcal D}$ be given. By part (b)(ii) of
Lemma \ref{LL} we may select a sequence $\{f_n\} \subset A_0$ so that
$\pi_\omega(d)\Omega_\omega = \lim_n\pi_\omega(f_n)\Omega_\omega.$
But then $$\omega(d) =
\langle\pi_\omega(d) \Omega_\omega,
\eta \rangle =
\lim_n \, \langle\pi_\omega(f_n) \Omega_\omega,
\eta \rangle = \lim_n \, \omega(f_n) = 0$$as required.

The equivalence with (iv) is now obvious.
\end{proof}

\begin{theorem}
Let $A$ be a maximal subdiagonal algebra.
Then $A$ satisfies the F
\& M Riesz property if and only if $\mathrm{dim}({\mathcal D}) < \infty$.
\end{theorem}

\begin{proof}  We proved the one direction in Proposition
\ref{fmrfd}.  For the other, let $\omega$ be a bounded
linear functional on $M$ which
annihilates $A_0$, and let $\omega_n$ and $\omega_s$ be the normal
and singular parts of $\omega$.  Write  $\omega_n =
\tau(h\cdot)$, for some $h \in L^1(M)$. We extend $\omega$,
$\omega_n$, and $\omega_s$, uniquely to normal functionals on the
enveloping von Neumann algebra (the double commutant in the
universal representation) and define $|\omega|$, $|\omega_n|$, and
$|\omega_s|$, to be the absolute values of these extensions
restricted to $M$. Then from for example (\cite{Hen}, cf.
\cite[Proposition 7]{E2}) applied to $\omega$ and $\tau$, we have
that as functionals on $M$, $|\omega_n|$ and $|\omega_s|$ are
respectively the normal and singular parts of $|\omega|$, and that
$|\omega| = |\omega_n| + |\omega_s|$. We note from
\cite[p.\ 270]{Dix} that there is no danger of confusion as regards
the absolute value of $\omega_n$ since the absolute value of
$\omega_n$ as a functional on $M$ and as a functional on the
enveloping von Neumann algebra coincide on $M$.
Now consider the positive functional $\rho$ given by $$\rho =
\tau + |\omega|.$$

Let $(\pi_{\rho}, \mathfrak{h}_{\rho}, \Omega_{\rho})$ be the GNS
representation constructed from $\rho$, and define $\rho_0$
by $\rho_0(a) = \langle\pi_{\rho}(a)(\Omega_{\rho}-\Omega_0),
\Omega_{\rho}-\Omega_0 \rangle$, where $\Omega_0$ is the orthogonal
projection of $\Omega_{\rho}$ onto the closure of
$\{\pi_{\rho}(a)\Omega_{\rho} : a \in A_0\}$. For any $f \in A_0$
and any $d \in {\mathcal D}_+$, we have by construction that
\begin{eqnarray*}
\|\pi_{\rho}(d^{1/2})(\Omega_{\rho}-\pi_{\rho}(f)\Omega_{\rho})\|^2
 &=& \rho(|d^{1/2}(\I-f)|^2)\\
 &\geq& \tau(|d^{1/2}(\I-f)|^2)\\
 &=& \tau(d - df - f^*d + |d^{1/2}f|^2)\\
 &=& \tau(d + |d^{1/2}f|^2)\\
 &\geq& \tau(d).
\end{eqnarray*}
On selecting a sequence $\{f_n\} \subset A_0$ so that
$\pi_{\rho}(f)\Omega_{\rho} \mapsto \Omega_0$, it
follows that $\rho_0(d) = \Vert \pi_{\rho}(d^{1/2})(\Omega_{\rho}- \Omega_0) \Vert^2
\geq \tau(d)$.
Hence  $\rho_0$ is faithful on ${\mathcal D}$,
and $\Omega_{\rho} \neq \Omega_0$.  Thus we may apply all of
Lemma \ref{LL} to $(\pi_{\rho},
\mathfrak{h}_{\rho}, \Omega_{\rho})$.

Next notice that for each $a$ in the enveloping von Neumann
algebra we have $$|\omega(a)|^2 \leq \|\omega\||\omega|(a^* a) \leq
\|\omega\|\rho(a^* a).$$  Thus on restricting
to elements of $M$,
and employing a standard argument, this implies that there exists a
vector $\eta \in \mathfrak{h}_{\rho}$ such that
 $$\omega(\cdot) =
\langle\pi_\rho(\cdot) \Omega_\rho,
\eta \rangle.$$Now consider the related functional
$$\widetilde{\omega}(\cdot) =
\langle\pi_\rho(\cdot) (\Omega_\rho-\Omega_0), \eta \rangle.$$
Select a sequence $\{f_n\} \subset A_0$ so that
$\pi_\rho(f_n)\Omega_\rho \rightarrow \Omega_0$. Let $a \in A_0$
be given. Since $A_0$ is an ideal, and since
$\omega$ annihilates $A_0$, we conclude that
\begin{eqnarray*}
\widetilde{\omega}(a) &=&
\langle\pi_\rho(a) (\Omega_\rho-\Omega_0),
\eta \rangle\\
 &=& \lim_n \, \langle\pi_\rho(a(\I-f_n)) \Omega_\rho,
 \eta \rangle\\
 &=& \lim_n \, \omega(a(\I-f_n))\\
 &=& 0.
\end{eqnarray*}
Thus $\widetilde{\omega}$ also annihilates $A_0$.

By part (c) of the Lemma we can find a sequence $\{a_n\}
\subset A$ such that
$\pi_{\rho}(a_n)(\Omega_{\rho}-\Omega_0) \, \to \, p_0\Omega_{\rho}$.
Let $a \in A_0$
be given. Since $A_0$ is an ideal, and since
$\widetilde{\omega}$ annihilates
$A_0$, we may now conclude that
\begin{eqnarray*}
\omega_n(a) &=& \langle\pi_\rho(a) p_0\Omega_\rho,
 \eta \rangle\\
 &=& \lim_n\langle\pi_\rho(aa_n) (\Omega_\rho-\Omega_0),
 \eta \rangle\\
 &=& \lim_n\widetilde{\omega}(aa_n)\\
 &=& 0.
\end{eqnarray*}
Thus $\omega_n$ annihilates $A_0$. But then so does $\omega_s =
\omega - \omega_n$. It now follows from Corollary \ref{eqfm}
that $A$ satisfies the F \& M Riesz property.
\end{proof}

\begin{corollary}   \label{hfmco}  If
$A$ is a maximal subdiagonal algebra with ${\mathcal D}$ finite
dimensional, and if $\omega \in M^*$ annihilates $A + A^*$, then
$\omega$ is singular.
\end{corollary}

\begin{proof}
Since $A$ satisfies the F \& M Riesz property, $\omega_n$ annihilates $A$.
Similarly, since $A^*$ satisfies the F \& M Riesz property, $\omega_n$ annihilates $A^*$.
Since $A$ is  subdiagonal, $\omega_n = 0$.
\end{proof}

\begin{corollary}   \label{ifF}  If $A$ has the F \& M Riesz property,
then any positive functional on $M$ which annihilates $A_0$ is normal.
\end{corollary}

\begin{proof}
If $\omega$ is a state on $M$ which annihilates $A_0$,
and if $A$ has the F \& M Riesz property, then the (positive) singular
part of $\omega$ is $0$ since it must annihilate $\I$.
\end{proof}



\section{The Gleason-Whitney theorem}


We say that an extension in $M^\star$ of a functional  in $A^\star$ is
a {\em Hahn-Banach extension} if it has the same norm.
If $A$ is a weak* closed subalgebra of $M$
then we say that  $A$ has property (GW1) if every
Hahn-Banach extension to $M$ of any normal functional on $A$,
is normal on $M$.
We say that  $A$ has property (GW2)  if there is at most one normal
Hahn-Banach extension to $M$ of
any normal functional on $A$.
We say that  $A$ has the {\em Gleason-Whitney property} (GW)
if it possesses (GW1) and (GW2).  This is simply saying
that there is a unique Hahn-Banach extension to $M$
of any normal functional on $A$, and this extension is normal.
Of course normal functionals on $A$ or on $M$ have to be
of the form $\tau(g \, \cdot \,)$ for some $g \in L^1(M)$.

\begin{theorem}  \label{GW}  If $A$ is a tracial
subalgebra of $M$ then $A$ is maximal subdiagonal
if and only if it possesses property (GW2).
If ${\mathcal D}$ is finite dimensional, then
$A$ is maximal subdiagonal
if and only if it possesses property (GW).
\end{theorem}

\begin{proof}    Suppose that $A$ possesses property (GW2).
To show that $A$ is maximal subdiagonal,
it suffices to show that if  $g \in L^1(M)$, with $\tau(g (A + A^*)) = 0$,
then $g = 0$.  By considering real and imaginary parts we may assume that
$g = g^*$.  Then $\tau(|g| \cdot)$ and $\tau((|g| + g) \cdot)$ are
positive normal functionals on $M$ which agree on $A$.
They are also Hahn-Banach extensions,
since the norm of a positive  functional is achieved at $1$.
Thus by (GW2), these functionals agree on $M$, and so $|g| + g = |g|$.
That is,  $g = 0$.

In the remainder of the proof
 suppose that $A$  is maximal subdiagonal.
Suppose that
$f, g \in L^1(M)$ correspond to two normal
Hahn-Banach extensions to $M$ of a given functional
on $A$.  Then
$\Vert f \Vert_1 = \Vert g \Vert_1$, and this quantity
equals the norm of the restriction to $A$.
We have $\tau((f-g)A) = 0$; since $A$ is subdiagonal
it follows from \cite[Lemma 4]{Sai} that $h = g - f \in
[A_0]_1$.   In order to establish (GW2), we need to show that $h = 0$.
Since Ball$(A)$ is weak* compact,
and since  $\Vert f \Vert_1$ equals
the norm of the above-mentioned restriction to $A$,
there exists
$a \in A$ of norm $1$ with $\tau(f a) = \Vert f \Vert_1$.
It is evident that
$$|a f|^2 = f^* a^* a f \leq f^* f = |f|^2 .$$
Now $0 \leq T \leq S$ in $L^p(M)$ implies
that $T^{\frac{1}{2}} \leq S^{\frac{1}{2}}$ (see e.g.\
\cite[Lemma 2.3]{Sc}, and we thank David Sherman for this
reference).  It follows that
$|a f| \leq |f|$.  On the other hand,
$\tau(|f|) = \tau(f a) = \tau(a f) \leq  \tau(|a f|)$.
Thus $\Vert |f| - |af| \Vert_1 = \tau(|f| - |a f|) = 0$, and so $|f| = |a f|$.
The functional $\psi = \tau(af \cdot)$ on $M$ must be positive
since $\psi(\I) = \tau(af) = \tau(|f|) = \tau(|af|) = \Vert \psi \Vert$.
Thus $af \geq 0$, and $af = |af| = |f|$.

Since $h \in [A_0]_1$ we have
$$\tau((f + h) a) = \tau(f a) = \Vert f \Vert_1 = \Vert g \Vert_1 = \Vert f + h
\Vert_1 .$$
An argument similar to that of the last paragraph shows that
$a (f + h) = |f + h|  \geq 0$.  Thus $a h$ is self-adjoint.   Since
$h \in [A_0]_1$ it is easy to see that $\tau(a h A) = 0$. Therefore from the 
self-adjointness of $ah$ one may deduce that
$\tau(a h (A + A^*)) = 0$.   Because $A$  is  subdiagonal, it follows
that $a h = 0$.
Thus
$$|f| = af = a (f+h) = |f + h| .$$
Let $e$ be the left support projection of $a$.  Then $e^\perp$ is the
projection onto Ker$(a^*)$.  We have
$|f| e^\perp = f^* a^* e^\perp = 0$. It follows that $f e^\perp = 0$.  Thus
$$0 = e^\perp f^* f e^\perp = e^\perp |f + h|^2 e^\perp
= e^\perp (f+h)^* (f+h) e^\perp = e^\perp h^* h e^\perp .$$
Hence $h e^\perp = 0$.
 To show that $h e = 0$, we reproduce the ideas
in the argument in the second paragraph of the proof.   Namely, note that
$|(f a)^*|^2 \leq |f^*|^2$ , so that
$|(f a)^*| \leq |f^*|$.
But $\tau(|f^*|) = \Vert f \Vert_1
= \tau(f a) \leq \tau(|(f a)^*|)$,
and as before this shows that $|(f a)^*| = |f^*|$.
Then also $\tau(f a) = \tau(|(f a)^*|)$, and as before this
shows that $f a  \geq 0$.  Similarly, $(f  + h) a \geq 0$.
So $h a$ is again selfadjoint, and this implies as before that $h a = 0$.
Thus $h  e = 0$, and so $h = he + he^\perp = 0$ as
required.

Now suppose that, in addition,  ${\mathcal D}$ is finite dimensional,
and that $\rho$ is a Hahn-Banach extension of a
normal functional $\omega$
on $A$.  By basic functional analysis,
 $\omega$ is the restriction of  a
normal functional $\tilde{\omega}$ on $M$.
  We may write $\rho = \rho_n + \rho_s$,
where $\rho_n$ and $\rho_s$ are respectively the normal and singular
parts, and $\Vert \rho \Vert = \Vert \rho_n  \Vert + \Vert \rho_s \Vert$.
Then $\rho - \tilde{\omega}$
annihilates $A$, and hence by our F.\ and M.\ Riesz theorem
both the
normal and singular
parts, $\rho_n - \tilde{\omega}$
and $\rho_s$ respectively,
annihilate $A_0$.  Hence they annihilate $A$, and in particular
$\rho_n = \omega$ on $A$.  But this implies that
$$\Vert \rho_n \Vert
+ \Vert \rho_s \Vert = \Vert \rho \Vert = \Vert \omega  \Vert
\leq \Vert \rho_n \Vert .$$
We conclude that $\rho_s  = 0$.   Thus $A$ also satisfies (GW1), and
hence (GW).
\end{proof}

\medskip

There is another (simpler) variant of the Gleason-Whitney theorem
\cite[p.\ 305]{Ho}, which transfers more easily to our setting:

\begin{theorem}  \label{GW2}  Let  $A$ be a maximal subdiagonal subalgebra of $M$
with ${\mathcal D}$ finite dimensional.
If $\omega$ is a normal functional on $M$ then $\omega$ is the unique
Hahn-Banach extension of its restriction to $A + A^*$.  In particular,
$\Vert \omega \Vert = \Vert \omega_{\vert A + A^*} \Vert$ for any $\omega \in M_*$.
\end{theorem}

\begin{proof}  Let $\rho$ be a Hahn-Banach extension of
the restriction of  $\omega$ to $A + A^*$.
We may write $\rho = \rho_n + \rho_s$,
where $\rho_n$ and $\rho_s$ are respectively the normal and singular
parts, and $\Vert \rho \Vert = \Vert \rho_n  \Vert + \Vert \rho_s \Vert$.
Then $\rho - \omega$ annihilates $A + A^*$.
By Corollary \ref{hfmco},
$\rho_n - \omega = (\rho  - \omega)_n = 0$.  As in the last part of the
previous proof, this implies that $\rho_s  = 0$.  So $\rho = \rho_n = \omega$.
\end{proof}

\medskip

{\bf Remark.}  If $g \in L^1(M)$, and  $\omega = \tau(g \, \cdot \, )$,
then the last result shows that 
$\Vert g \Vert_1$ is the norm of the restriction of $\omega$ to $A + A^*$.

\begin{corollary}  \label{Kap} {\rm
(Kaplansky density theorem for subdiagonal algebras) } \
 Let  $A$ be a maximal subdiagonal subalgebra of $M$
with ${\mathcal D}$ finite dimensional.
Then the unit ball of $A + A^*$ is weak* dense in ${\rm Ball}(M)$.
\end{corollary}

\begin{proof}  If $C$ is the unit ball of $A + A^*$, it follows from the
last remark that the pre-polar of $C$ is ${\rm Ball}(M_\star)$.
By the bipolar theorem, $C$ is  weak* dense in ${\rm Ball}(M)$.
\end{proof}

\medskip

{\bf Remark.}  We do not know if the last few results hold
without the assumption that ${\mathcal D}$  be finite dimensional.

\section{Szeg\"o and Kolmogorov theorems for $L^p(M)$}

Arveson formulated the Szeg\"o theorem for $L^2(M)$ in terms
of the Kadison-Fuglede determinant $\Delta(\cdot)$.  The long-outstanding
open question of whether general maximal
subdiagonal algebras satisfy the Szeg\"o theorem for $L^2(M)$,
was eventually settled in the affirmative in \cite{Lab}.
We will now extend this result to $L^p(M)$.  We refer the reader
to \cite{AIOA,BL2} for the properties of the 
Kadison-Fuglede determinant which we shall need.

\begin{lemma} \label{Szpl}  $\Delta(b^p) = \Delta(b)^p$ for $p \geq 1$ and
$b \in M_+$.  \end{lemma}

\begin{proof}  By the multiplicativity property of $\Delta$, the
relation clearly holds for dyadic rationals.
 We may assume that $0 \leq b \leq 1$.
In this case, by the functional calculus it is clear that
$b^q \leq b^p$ if $0 < p \leq q$.  If $q$ is any dyadic rational
bigger than $p$ then $$\Delta(b)^q = \Delta(b^q) \leq \Delta(b^p) .$$
It follows that  $\Delta(b)^p \leq \Delta(b^p)$.
Replacing $p$ by $1/p$, we have
$\Delta(b^p)^{\frac{1}{p}} \leq \Delta((b^p)^{\frac{1}{p}}) = \Delta(b)$,
which gives the other direction.   \end{proof}

\begin{theorem} \label{Szp} {\rm (Szeg\"o theorem for $L^p(M)$)}
 \  Suppose that $A$ is maximal  subdiagonal,
and  $1 \leq p < \infty$.
 If $h \in L^1(M)_+$ then $\Delta(h) = \inf \{ \tau(h |a+d|^p)
: a \in A_0, d \in {\mathcal D} , \Delta(d) \geq 1 \}$.
\end{theorem}

\begin{proof}   We set
$${\mathcal S}_p = \{ |a|^p :  a \in A , \,
\Delta(\Phi(a)) \geq 1 \} , $$
$$ {\mathcal S} = \{ a^* a : a \in A^{-1} , \,
\Delta(a) \geq 1 \} .$$
By the modification in \cite[Proposition 3.5]{BL2}
of a trick of Arveson's from
  \cite[Theorem 4.4.3]{AIOA},
it suffices to show that
the closure of ${\mathcal S}_p$ equals the closure of
${\mathcal S}$.
First we show that ${\mathcal S} \subset {\mathcal S}_p$.  Indeed,
if $b \in {\mathcal S}$ then $b$ is invertible, and therefore so
is $b^{\frac{1}{p}}$.  Since $A$ has factorization,
there is an $a \in A^{-1}$ with $|a| = b^{\frac{1}{p}}$.
By Lemma \ref{Szpl} and Jensen's formula \cite{AIOA,Lab}
we have $$\Delta(\Phi(a)) = \Delta(a)  = \Delta(|a|)
= \Delta(b^{\frac{1}{p}}) = \Delta(b)^{\frac{1}{p}} \geq 1 .$$
Hence $b = |a|^p \in {\mathcal S}_p$.

Suppose that  $b \in {\mathcal S}_p$.
If $b = |a|^p$ where $\Delta(\Phi(a)) \geq 1$ then
by Jensen's inequality \cite{AIOA,Lab}
we have $\Delta(a) = \Delta(|a|) \geq 1$.
Hence by Lemma \ref{Szpl} we have
 $\Delta(b) \geq 1$.  If $n \in \N$ then since $A$  has factorization,
there exists a  $c \in  A^{-1}$ with
$b + \frac{1}{n} 1 =  c^* c$.
 Thus $$\Delta(c)^2 =  \Delta(b + \frac{1}{n} 1) \geq  \Delta(b) \geq 1 .$$
Thus $b + \frac{1}{n} 1 = c^* c  \in {\mathcal S}$, and
we deduce that $b \in \overline{{\mathcal S}}$.
Hence
$\overline{{\mathcal S}_p} \subset \overline{{\mathcal S}}$.
\end{proof}

\medskip

Note that the following generalized Kolmogorov theorem
is not true for all maximal  subdiagonal algebras.
For example, take $A = M = L^\infty[0,1]$.

\begin{theorem} \label{Kolm}
Suppose that $A$ is an antisymmetric
maximal  subdiagonal algebra.
If $h \in L^1(M)_+$ then $\inf \{ \tau(h |\I + f|^2)
: f \in A_0 + A_0^* \}$ is either $\tau(h^{-1})^{-\frac{1}{2}}$,
if $h^{-1}$ exists in the sense of unbounded operators and is in
$L^1(M)$; or the infimum is $0$ if $h^{-1} \notin L^1(M)$.
More generally, if $1 \leq p < \infty$ then
 $\inf \{ \tau(|(\I+f) h^{\frac{1}{p}}|^p)
: f \in A_0 + A_0^* \}$ is either $0$ if $h^{-1} \notin L^{1/(p-1)}(M)$,
or $\tau(h^{-\frac{1}{p-1}})^{\frac{1}{p} - 1}$ if $h^{-1} \in L^{1/(p-1)}(M)$.
 \end{theorem}

\begin{proof}   We formally follow the proof of Forelli as adapted
in \cite[p.\ 247]{SW}.  Let $h \in L^1(M)_+$, and $1/p + 1/q = 1$.
Define $L^p(M,h)$ to be the completion
in $L^p(M)$ of $M h^{\frac{1}{p}}$.   Note that if $e$ is the
support projection of a positive $x \in L^p(M)$ then
it is well known (see e.g.\ \cite[Lemma 2.2]{JS}) that
 $L^p(M) e$ equals the closure in $L^p(M)$ of $M x$.
Hence $L^p(M,h) = L^p(M) e$, where $e$ is the support projection of  $h$.
Now for any projection $e \in M$ it is an easy exercise
 to prove that
the dual of $L^p(M) e$ is
 $e L^q(M)$ (see e.g.\ \cite{JS}).
It follows that  the dual of $L^p(M,h)$ is $L^q(M,h)$.

If $k \in L^p(M,h)$ then $k h^{\frac{1}{q}} \in L^p(M) L^q(M) \subset
L^1(M)$.   We view $A_0 + A_0^*$ in $L^p(M,h)$ as its image
$(A_0 + A_0^*) h^{\frac{1}{p}}$,
and let $N$ be the annihilator of this in $L^q(M,h)$.  That is, $g \in
N$ iff $g \in L^q(M,h)$ and $$0 =
\tau(h^{\frac{1}{p}} (A_0 + A_0^*) g) =
\tau((A_0 + A_0^*) g h^{\frac{1}{p}}) .$$
Since $gh^{\frac{1}{p}} \in L^1(M)$
the last equation holds iff $gh^{\frac{1}{p}} = c \I$,
where $c$ is  a constant.    Since $h$ is selfadjoint, if $c \neq 0$ then
it follows that
$h^{-\frac{1}{p}}$ exists in the sense of unbounded operators, and its
closure is the constant multiple $d g \in L^q(M)$,
where $d = c^{-1}$.  (Since we are in the finite case,
there is no difficulty with $\tau$-measurability here,
this is automatic \cite{Terp}.)
 If $c = 0$ then
$gh^{\frac{1}{p}} = 0$ which implies that $g = 0$.  To see the
last statement note that if $h^{\frac{1}{p}}$ is viewed as
a selfadjoint unbounded operator on a Hilbert space $H$, and if
$e$ is its support projection, which equals the support projection
of $h^{\frac{1}{q}}$, then $e h^{\frac{1}{p}} = h^{\frac{1}{p}}$,
and so $h^{\frac{1}{p}} e = h^{\frac{1}{p}}$.  Since
$g \in \overline{M h^{\frac{1}{q}}}$, we have $g e = g$.  However $g e
 = 0$ since $g h^{\frac{1}{p}} = 0$.
Thus if $g$ has norm 1 then $c \neq 0, h^{-\frac{1}{p}} \in L^q(M)$ and
$|d| = \Vert h^{-\frac{1}{p}} \Vert_{L^q(M)} = \tau(h^{-\frac{q}{p}})^{\frac{1}{q}}
$.

The infimum in the theorem is the $p$th power of the norm of $\I$ in the
quotient space of $L^p(M,h)$ modulo the closure of $A_0 + A_0^*$.  Since the
dual of this quotient is $(A_0 + A_0^*)^\perp = N$, this infimum
equals the $p$th power of
$\sup \{ | \tau(g h^{\frac{1}{p}}) | : g \in N , \Vert g
\Vert_{L^q(M)} \leq 1 \} $.  This equals $0$ if no $g \in N$ has norm 1; otherwise it
equals $\tau(h^{-\frac{q}{p}})^{-\frac{1}{q}}
= \tau(h^{-\frac{1}{p-1}})^{-\frac{1}{q}}$
by the above.
Indeed, the infimum is $0$ iff $\tau(gh^{\frac{1}{p}}) = 0$ for all $g \in N$.
Since $gh^{\frac{1}{p}}$ is constant, this occurs iff
$gh^{\frac{1}{p}} = 0$, which as we saw above happens iff
$g = 0$.  Thus the infimum is $0$ iff $N = (0)$ iff
$(A_0 + A_0^*) h^{\frac{1}{p}}$ is dense in
$L^p(M,h)$.  Since $h^{\frac{1}{p}} \in L^p(M,h)$, the
latter condition implies that there is a sequence $(g_n)$ in
$A_0 + A_0^*$ with $g_n h^{\frac{1}{p}} \to h^{\frac{1}{p}}$ in
$p$-norm.  If $h^{-1/p} \in L^q(M)$ then by H\"older's inequality
we have $\tau(|g_n - \I|) \to 0$, which is impossible
since $1 = |\tau(g_n - \I)| \leq \tau(|g_n - \I|)$.   \end{proof}

\medskip


{\bf Acknowledgements.}  We thank Mike Marsalli for many valuable discussions,
and Marius Junge for a helpful insight concerning Theorem \ref{ex}.

\medskip





\end{document}